\documentclass[12pt]{amsart}

\usepackage{latexsym} 
\usepackage{amsthm}
\usepackage{amsfonts} 
\usepackage{graphicx}
\usepackage[abs]{overpic}

\newtheorem{theorem}{Theorem}[section]
\newtheorem{lemma}[theorem]{Lemma}
\newtheorem{proposition}[theorem]{Proposition}
\newtheorem{corollary}[theorem]{Corollary}

\theoremstyle{definition} 
\newtheorem{definition}[theorem]{Definition}

\theoremstyle{remark}

\numberwithin{equation}{section}

\begin{document}

\markboth{Kenneth C. Millett, Michael Piatek, and Eric J. Rawdon}
{Polygonal knot space near ropelength-minimized knots}

\title{Polygonal knot space near ropelength-minimized knots}

\author{Kenneth C. Millett} 
\address{Department of Mathematics,
University of California, Santa Barbara, Santa Barbara, CA 93106, USA,
email: millett@math.ucsb.edu} 

\author{Michael Piatek}
\address{Department
of Computer Science and Engineering, University of Washington, Seattle,
WA 98195, USA, email: piatek@cs.washington.edu} 

\author{Eric Rawdon} 
\address{Department
of Mathematics and Computer Science, Duquesne University, Pittsburgh,
PA 15282, USA, email: rawdon@mathcs.duq.edu} 

\maketitle

\begin{abstract}
For a polygonal knot $K$, it is shown that a tube of radius $R(K)$,
the polygonal thickness radius, is an embedded torus.  Given a thick
configuration $K$, perturbations of size $r<R(K)$ define satellite
structures, or local knotting.  We explore knotting within these tubes
both theoretically and numerically.  We provide bounds on perturbation
radii for which we can see small trefoil and figure-eight summands and
use Monte Carlo simulations to approximate the relative probabilities
of these structures as a function of the number of edges.
\end{abstract}



\section{Introduction}
Equilateral polygonal knots have been used widely as statistical
models of macromolecules such as DNA.  For a fixed number of edges $n$
these knots inhabit a space $Equ(n)$ whose structure reflects the
possibility of changes in conformation subject to constraints such as
the preservation of edge lengths and the prohibition of singularities
such as one edge passing through another.  This space is a subspace of
the space of polygonal knots $Geo(n)$ in which the edges may change in
length but not in number and again, edges may not pass through each
other.  In this paper, we will study the behavior of equilateral knots
under spatial perturbation of their vertices in $Geo(n)$.  These
perturbations may be thought of as changes in the configuration during
which we suspend the constraint of edge length preservation and, when
the scale of perturbation is sufficiently large, the prohibition
against edges passing through other edges.  We will require that the
number of edges remains the same.  Of special interest are the
``thick'' equilateral knots, representatives of the knot type in which
the edges are quite separated.  We explore the structure of local
knotting near these ``thick'' knots which, by our definition, will be
those which occur from perturbations within a thick tube surrounding
the knot.  Global knotting changes can occur when the size of the
perturbation exceeds the radius of this tube.

Key to this study is the ``Tube Theorem'' (Theorem \ref{tubethm} in
Section \ref{tubethmsection}) describing the explicit structure of a
uniformly large embedded solid torus containing the knot as its
center.  The radius of this tube about the knot is the thickness
radius of the knot conformation $R(K)$ defined in \cite{meideal,mine}.
The theorem and the associated value of $R(K)$ are the starting point
of our analysis.  The precise definition of $R(K)$ and a discussion of
some of its properties are contained in the next section.

Perturbations of the knot that remain within this tube define
satellites of the knot.  The nature of these perturbations vary from
being equivalent polygonal knots, for very small perturbations, to
those which are connected sums with small trefoils (a ``local knot'')
and, finally, to those satellites that may constitute a more complex
structure distributed across the entire tube but whose structure still
derives from a local perturbation of the core knot.  Once the scale of
the perturbation allows the edges to pass beyond the boundary of the
tube, crossings between arclength-distant edges in the knot may occur
and a ``global'' change in knotting is possible, but for the purpose
of this note, we concentrate on the perturbations causing local
knotting changes.

From the perspective of $Equ(n)$, as a subspace of $Geo(n)$, the
perturbed knots can be understood as nearby neighbors of the knot $K$
in knot space and can be expected to have structure that bears a close
resemblance to that of $K$.  We will see that the critical
perturbation parameters that quantify the extent of structural change
are the thickness radius of the knot and the number and length of its
edges.  So long as the scale of the perturbation is sufficiently small
in terms of these quantities, the resulting knot will be geometrically
equivalent, that is lying in the same connected component of $Geo(n)$
as does the knot $K$.  If the relationship between the number of
edges, edge length, and the thickness radius allows, the neighbors may
contain small trefoil summands or summands of other small knots.  This
relationship forms the basis of the passage from a ``rigid'' knot
regime to one whose behavior more closely resembles that of smooth
knots in a thickest state.  That is to say, for some equilateral
knots, it is impossible to form local knots even in its thickest
conformation due to the small number of edges.  The minimal number of
edges required to form a local knot under perturbation in an
equilateral knot type is one possible indicator of the passage from
being a ``rigid'' knot to a structure that permits sufficient
flexibility of motion to better model the behavior of large knots.

\begin{figure}
  \centering
  \includegraphics[height=1.2in]{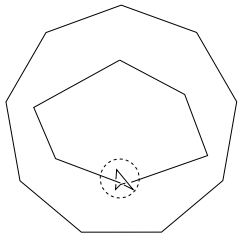}
  \vspace*{8pt}
  \caption{A nine edge trefoil satellite of the trivial knot.}
  \label{fig1}
\end{figure}

\begin{figure}
  \centering
  \includegraphics[height=1.2in,angle=90]{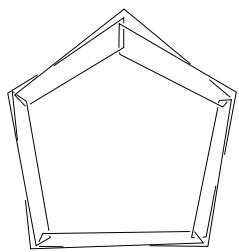}
  \vspace*{8pt}
  \caption{A complex 25 edge satellite of the trivial knot.}
  \label{fig2}
\end{figure}

For sufficiently large knots, as expressed in terms of the
relationship between the number of edges and the knot's thickness
radius, it is possible to encounter perturbations that are organized
across its entire length so as to produce more complex satellite
structures.  Such perturbations require more edges than are necessary
to produce local knotting and provide a second point at which the
character of the knot may be thought of as changing.  For example, for
the unknot, nine edges are required in order that a (rather large)
perturbation within the tube gives rise to a local trefoil as shown in
Figure \ref{fig1}.  By way of contrast, a rough analysis estimates
that at least 25 edges are required for a satellite that is not a sum
of local knots as shown in Figure \ref{fig2}.  So long as the
relationship between the thickness radius and the edge lengths permit,
we will see that the perturbations of edges can remain within the tube
and exhibit increasingly complex local knotting structure and an
exponential growth in occurrence and complexity as expected. This
structure is characteristic of the neighborhood, as measured by the
perturbation scale, of the $K$ in knot space.

In this paper we will provide both a theoretical analysis and Monte
Carlo simulation data that illuminate the local structure of knot
space as measured by the perturbation scale associated to an
equilateral knot.  The theoretical dimension will provide estimates of
the scale of perturbation required to attain different satellite
structures while the Monte Carlo data will relate these estimates to
observation.  Using a random sampling of the set of perturbations and
a HOMFLY polynomial analysis of this data, we are able to study and
characterize the complexity of knots nearby the knot determining core
of the tube.  One may think of this effort as an indirect strategy to
study the local geometry and topology of knot space.  For example, we
analyze the number of edges required to obtain different satellite
structures (usually occurring as connected sums of prime knots) as
well as the growth in the relative frequency as a function of the
number of edges.  These are exhibited in the differences in the data
at the local scale near the thickest trivial, trefoil, and
figure-eight knots.

\section{Background on Ropelength}

We explore perturbations about ropelength-minimized equilateral
polygonal knots, the so-called \textit{ideal knots}.  For
completeness, we will include background on ropelength for smooth
knots, which for the sake of this discussion means $C^2$, and for
polygonal knots, in which edge lengths are no longer required to be
equal.

The ropelength of a knot was defined in \cite{bosimple} and the basic
theory developed in \cite{lsdr}, although similar ideas appeared
previously \cite{federer,otto,nab}.  Smooth ropelength models an
idealize rope as a non self-intersecting tube with a knot as the
filament-like core.

\begin{definition}
For a $C^{2}$ smooth knot $K$ and $x\in K$, let $D_r(x)$ be the disk
of radius $r$ centered at $x$ lying in the plane normal to the tangent
vector at $x$.  Let
\begin{equation*}
R(K)=\sup\{r>0\,:\,D_r(x)\cap D_r(y)=\emptyset
\text{ for all }x\not = y \in K\}.
\end{equation*}
The quantity $R(K)$ is called the {\it thickness radius} of $K$.
Define the {\it ropelength} of
$K$ to be 
\begin{equation*}
Rope(K)=Length(K)/R(K)
\end{equation*}
where $Length(K)$ is the
arclength of $K$.
\end{definition}

For a fixed knot configuration, $R(K)$ is the radius of a thickest
tube that can be placed about $K$ without self-intersection.  Given a
knot surrounded by an impenetrable tube of some radius, there are two
types of interactions between normal disks that restrict the
conformations that can be obtained.  First, the tube cannot bend too
quickly, a restriction on the curving of the core.  Second, two
arclength-distant points in the core cannot be any closer than twice
the radius of the tube, a restriction on the distance between pairs of
points bounded away from the diagonal of $K\times K$.  These intuitive
observations are captured by the quantities below and the subsequent
lemma.

\begin{definition}
For a $C^2$ knot $K$ with unit tangent map $T$, let $MinRad(K)$ be the
minimum radius of curvature.  The {\it doubly-critical self-distance}
is the minimum distance between pairs of distinct points on the knot
whose chord is perpendicular to the tangent vectors at both of the
points.  In other words, let
\begin{equation*}
DC(K)=\{(x,y)\in K\times K\,:\,T(x)\perp\overline{xy}\perp T(y), 
x\not = y\}
\end{equation*}
where $\overline{xy}$ is the chord connecting $x$ and $y$.
Define the {\it doubly-critical self-distance} by
\begin{equation*}
dcsd(K)=\min\{\Vert x-y\Vert\,:\,(x,y)\in DC(K)\},
\end{equation*}
where $\Vert \cdot \Vert$ is the standard $\mathbb{R}^3$ norm.
\end{definition}

There is a fundamental relationship between $R(K)$, $MinRad(K)$, and
$dcsd(K)$.

\begin{lemma}
Suppose $K$ is a $C^2$ knot.  Then
$$R(K)=\min\left\{MinRad(K),dcsd(K)/2\right\}.$$
\label{thicknessc2}
\end{lemma}
\begin{proof}
See \cite{lsdr}.
\end{proof}

A surprising result from \cite{lsdr} establishes that, in fact, the
doubly-critical self-distance can be replaced by the singly-critical
self-distance, an idea attributed to J.~O'Hara and N.~Kuiper
\cite{oharafamily}.  We prove a similar result for polygons in the
next section.

\begin{proposition}
Let $$SC(K) = \{(x,y)\in K\times K\,:\, T(x)\perp \overline{xy}, x\not
= y\}$$ and the {\it singly-critical self-distance} be 
$scsd(K)=\min_{(x,y)\in SC(K)}\Vert x-y\Vert$.  Then
$$R(K) = \min\{MinRad(K),scsd(K)/2\}\,.$$
\label{scsdthickness}
\end{proposition}
\begin{proof}
See \cite{lsdr}.
\end{proof}

The thickness radius and ropelength definitions can be extended to
$C^{1,1}$ curves (see \cite{maddocks,cksminimum,lsdr2}).  Since $R(K)$
changes with scale and we are mainly interested in the ``shape'' of
the optima, the scale-invariant ropelength is the quantity typically
studied.

In \cite{meideal,mine}, the polygonal thickness radius and polygonal
ropelength functions were defined in the spirit of the
characterization in Lemma \ref{thicknessc2}.  Before we define this,
we introduce some notation.

\begin{figure}
  \begin{center}
    \begin{overpic}{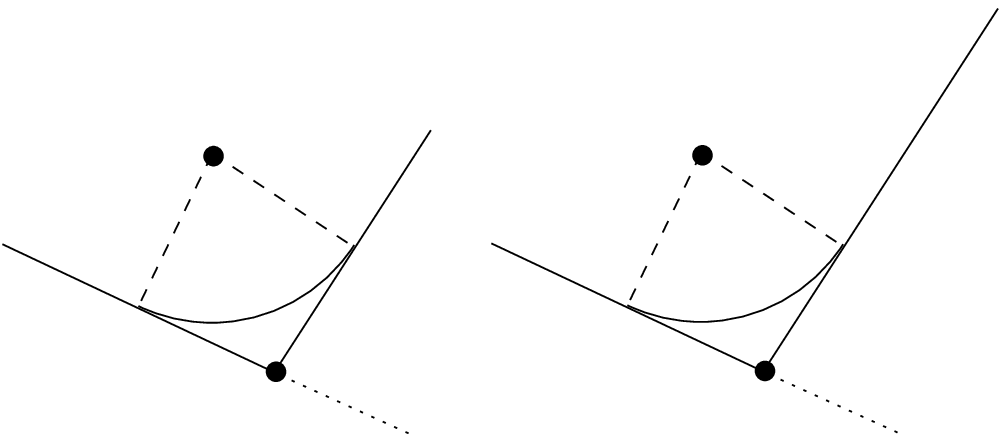}
      \put(90,20){$angle(v)$}
      \put(230,20){$angle(v)$}
      \put(10,58){$Rad(v)$}
      \put(152,58){$Rad(v)$}
      \put(75,7){$v$}
      \put(215,7){$v$}
    \end{overpic}
  \end{center}
  \caption{An arc of a circle of radius $Rad(v)$ can be inscribed so
    that the arc is tangent at the midpoint of the shorter adjacent
    edge.  On the left, the two adjacent edges have identical length,
    so the arc intersects both edges tangentially at the midpoints.
    On the right, the arc intersects the longer edge at a point short
    of the midpoint.}
\label{minradpic}
\end{figure}

For an $n$-edge polygonal knot $K$, define the following:
\begin{itemize}
\item Let $v_1,\cdots,v_{n}$ be the vertices of $K$.
For convenience, we implicitly take all subscripts modulo $n$.
\item Let $e_1,\cdots,e_{n}$ be the edges of $K$, where $e_i$ is the
edge connecting $v_i$ to $v_{i+1}$.  
\item Let $|e_i|$ be the length of the edge $e_i$.
\item Let $angle(v_i)$ be the measure of the turning angle at $v_i$
(see Figure \ref{minradpic}).
\item Given a knot $K$ and $x\in K$, let $d_x:K\to \mathbb{R}$ be defined
by $d_x(y)=\Vert x-y\Vert$.
\end{itemize}

\begin{definition}
For a vertex $v_i$ on an $n$-edge polygonal knot $K$, let 
\begin{equation*}
Rad(v_i)= \frac{\min\{|e_{i-1}|,|e_i|\}}
{2\tan\left(angle(v_i)/2\right)}
\end{equation*}
and 
\begin{equation*}
MinRad(K)=\min\limits_{i=1,\cdots,n}Rad(v_i).
\end{equation*}
Note that $Rad(v_i)$ is the radius of a circular arc that can be
inscribed at $v_i$ so that the arc is tangent to both edges adjacent
to $v_i$ and the arc intersects the shorter adjacent edge at its
midpoint (see Figure \ref{minradpic}).  We call $y$ a {\it turning
point} for $x$ if $d_x$ changes from increasing to decreasing or from
decreasing to increasing at $y$.  Let
\begin{equation*}
DC(K)=\{(x,y)\in K\times K\,:\,x\not = y
\text{ turning points of }d_y \text{ and }d_x \text{ respectively}\}.
\end{equation*}
Define the {\it doubly-critical self-distance} of $K$ to be
\begin{equation*}
dcsd(K)=\min\{\Vert x-y\Vert\,:\,(x,y)\in DC(K)\}.
\end{equation*}
We call a pair $(x,y)\in DC(K)$ a {\it doubly critical pair}.
\label{polydefs}
\end{definition}

\begin{definition}
For a polygonal knot $K$, let 
\begin{equation*}
R(K)=\min\left\{
MinRad(K),dcsd(K)/2\right\}\text{ (polygonal thickness radius)}
\end{equation*}
and
\begin{equation*}
Rope(K)=\frac{Length(K)}{R(K)}\text{ (polygonal ropelength)}.
\end{equation*}
\end{definition}

In \cite{meideal,mine}, it was shown that, under mild geometric
hypotheses, for finer and finer inscribed polygonal approximations of
a smooth knot $K$, the thickness radius and ropelength of the polygons
converge to the respective values of $K$.  In \cite{mecancomputers},
we showed that the ropelength of polygonal optima converge to the
minimum smooth ropelength and that a subsequence of the minimizing
polygons converges pointwise to a smooth ropelength minimum.

\section{Self-distance and lemmas about polygonal arcs}

In this section, we prove that polygonal thickness radius can be
reformulated in terms of the singly-critical self-distance.  This is
the polygonal version of Proposition \ref{scsdthickness} and further
shows how closely the polygonal thickness radius models the behavior
in the smooth case.  We need some definitions and lemmas before we can
prove this result.

\begin{definition}
For a polygonal knot $K$, let $MD(K)$ be the 
minimum distance between non-adjacent edges of $K$ and
$MD(e,v)$ be the minimum distance between edge $e$ and
vertex $v$.
\end{definition}

Notice that $MD(K)$ is bounded above by the minimum edge length of
$K$, which we denote $MinEdge(K)$.

\begin{definition}
For a polygonal knot $K$, let $SC(K)$ consist of all pairs 
$(x,y)\in K\times K$ such that $x$ and $y$ are on non-adjacent edges
and $y$ is a turning point which is a local minimum of $d_x$.
Define the {\it singly-critical self-distance} of $K$ by
$$scsd(K) = \min_{(x,y)\in SC(K)}\Vert x-y\Vert\,.$$
\end{definition}

Since $SC(K)$ is a closed set, $scsd(K)$ is well-defined.  We need the
following three lemmas to recharacterize the polygonal thickness
radius with $scsd$ used in place of $dcsd$.  The first provides a
lower bound on the distance between an edge and the vertices of the
adjacent edges.

\begin{figure}
  \begin{center}
    \begin{overpic}[height=2in]{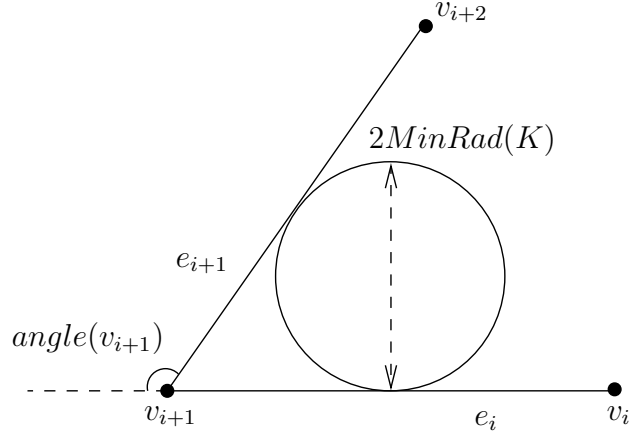}
      \put(-5,20){$angle(v_{i+1})$}
      \put(45,-8){$v_{i+1}$}
      \put(57,50){$e_{i+1}$}
      \put(130,95){$2MinRad(K)$}
      \put(155,145){$v_{i+2}$}
      \put(170,-10){$e_i$}
      \put(220,-8){$v_i$}
    \end{overpic}
  \end{center}
  \vspace*{10pt}
  \caption{When $angle(v_{i+1})\geq \pi/2$, the minimum distance from
    $v_{i+2}$ to $e_{i}$ is $\geq 2\,MinRad(K)$.}
  \label{onesep}
\end{figure}

\begin{lemma}
For a polygonal knot $K$, $$MD(e_i,v_{i+2}) \geq
\min\{2\,MinRad(K),MinEdge(K)\}.$$ In particular, when
$angle(v_{i+1})\geq \pi/2$, we have $$MD(e_i,v_{i+2}) \geq
2\,MinRad(K).$$ A similar result holds for $MD(e_i,v_{i-1})$.
\label{mu1sep}
\end{lemma}
\begin{proof}
If $angle(v_{i+1})\leq \pi/2$, then the minimum distance is realized at
$v_i$ and $v_{i+1}$ which is at least $MinEdge(K)$.  If
$angle(v_{i+1})>\pi/2$, then $\,2MinRad(K)\leq
\frac{|e_i|}{\tan\left(angle(v_{i+1})/2\right)}$ and
$d(v_i,e_{i+1})=|e_i|\,\sin(\pi-angle(v_{i+1}))$ (see Figure
\ref{onesep}).  Simple trigonometry show that
$2\,MinRad(K)\leq d(v_i,e_{i+1})$.
\end{proof}

The {\it total curvature} between two points $tc(x,y)$ on a polygon is
the minimum (over the two arcs connecting $x$ to $y$) of the sum of
the turning angles along the arc, including the angle at $x$ and/or
$y$ if one or both happen to be vertices.  The total curvature is used
to recharacterize polygonal $R(K)$ in Lemma \ref{polytcchar}.

\begin{lemma}
Suppose $x,y$ are points on a polygonal knot $K$ and
let $TC(K)=\{(x,y)\in K\times K\;:\;tc(x,y)\geq \pi\}$.
\begin{itemize}
\item[(a)] If $tc(x,y)<\pi$, then $(x,y)$ is not a doubly critical pair.
\item[(b)] If $tc(x,y)\geq \pi$, then $\Vert x-y\Vert \geq 2\,R(K)$.
\item[(c)] If $(x,y)\in TC(K)$ is on the boundary of $TC(K)$, then
$\Vert x-y\Vert \geq 2\,MinRad(K)$.
\end{itemize}
\label{tclemma}
\end{lemma}
\begin{proof}
See \cite{mine}.
\end{proof}

The following is a characterization of polygonal $R(K)$ that will be
used in the proof of Theorem \ref{polyscsd}, the polygonal version of
Proposition \ref{scsdthickness}.

\begin{lemma}
For a polygonal knot $K$,
\begin{equation*}
R(K)=\min\left\{MinRad(K),
\min_{(x,y)\in TC(K)}\Vert x-y\Vert/2\right\}.
\end{equation*}
\label{polytcchar}
\end{lemma}
\begin{proof}
See \cite{mine}.
\end{proof}

\begin{theorem}
For a polygonal knot $K$, $$R(K)=\min\left\{MinRad(K),
scsd(K)/2\right\}.$$
\label{polyscsd}
\end{theorem}
\begin{proof}

Usually, we will have $scsd(K) \leq dcsd(K)$.  However, recall that
$scsd(K)$ is defined only over pairs $(x,y)$ where $y$ is a local
minimum turning point of $d_x$ while $dcsd(K)$ is defined over all
doubly turning pairs, in particular those which are local maxima.

Suppose $scsd(K) > dcsd(K)$.  We show that in such a case $R(K) =
MinRad(K) \leq dcsd(K)/2 < scsd(K)/2$.  Suppose $dcsd(K)$ is realized
at the pair $(x_0,y_0)\in DC(K)$.  If $x_0$ is a local minimum for
$d_{y_0}$ or $y_0$ is a local minimum for $d_{x_0}$, then $scsd(K)
\leq dcsd(K)$.  Thus, $x_0$ and $y_0$ are local maxima for $d_{y_0}$
and $d_{x_0}$ respectively.  Note that $\min_{(x,y)\in TC(K)}\Vert
x-y\Vert \leq dcsd(K)$ by Lemma \ref{tclemma}(a).  If
$dcsd(K)/2<MinRad(K)$ then Lemma \ref{polytcchar} tells us that
$\min_{(x,y)\in TC(K)}\Vert x-y\Vert = dcsd(K)$.  In other words,
$(x,y)$ must be a global minimum of the distance function over
$TC(K)$.  By Lemma \ref{tclemma}(c), this global minimum cannot occur
on the boundary of $TC(K)$, and, thus, the points $x_0$ and $y_0$ must
be minima for $d_{y_0}$ and $d_{x_0}$ respectively. Thus,
$R(K)=MinRad(K)$ whenever $scsd(K) > dcsd(K)$, and the result holds.

So we can assume that $scsd(K)\leq dcsd(K)$.  Certainly, if
$MinRad(K)\leq scsd(K)/2\leq dcsd(K)/2$, then the result holds.  So
what remains is to show that if $scsd(K)/2<MinRad(K)$, then
$scsd(K)=dcsd(K)$.

Suppose $scsd(K)/2<MinRad(K)$ and $(x,y)$ is a pair in $SC(K)$ that
realizes $scsd(K)$.  In the following arguments, we assume that
$scsd(K)<dcsd(K)$ and obtain the contradiction that $scsd(K)<\Vert
x-y\Vert$.

Consider $B$, the ball of radius $scsd(K)$ centered at $y$.  We
consider the cases where $x$ and $y$ are vertex or non-vertex points.
When $x$ is a non-vertex point (i.e.~lies in the interior of an edge),
let $\alpha_x$ be the edge containing $x$ and define $\alpha_y$
analogously.

\noindent{\bf Case 1:} $y$, $x$ both non-vertex points\\ 
Since $(x,y)\in SC(K)$, $\alpha_x$ and $\alpha_y$ are non-adjacent.
Since $(x,y)$ is not doubly critical, $\alpha_x$ pierces $B$ at $x$.
Since $y$ is interior to $\alpha_y$, $y$ divides $\alpha_y$ into two
pieces, say $\alpha^1_y$ and $\alpha_y^2$.  Since $(x,y)\in SC(K)$,
$\overline{xy}\perp \alpha_y$ at $y$, so $\Vert x-y\Vert $ is a local
minimum of $d_x$.  Thus, there exist $y_1\in \alpha_y^1$ and
$y_2\in\alpha_y^2$, both interior to $\alpha_y$ such that
$$\Vert x-y_1\Vert >\Vert x-y\Vert \text{ and } \Vert x-y_2\Vert
>\Vert x-y\Vert \,.$$ 
Furthermore, there exists an $x'$ arbitrarily
close to $x$ on $\alpha_x\cap B$ such that
$$\Vert x'-y\Vert <\Vert x-y\Vert , \Vert x'-y_1\Vert >\Vert x-y\Vert , \text{ and } \Vert x'-y_2\Vert >\Vert x-y\Vert .$$
Thus, there exists a point $y'$ between $y_1$ and $y_2$ (interior
to $\alpha_y$) that minimizes the distance to $x'$.  Then
$(x',y')\in SC(K)$ with $\Vert x'-y'\Vert <\Vert x-y\Vert =scsd(K)$,
which is a contradiction.

\noindent{\bf Case 2:} $y$ is a vertex and $x$ is non-vertex\\
Then $\alpha_x$ pierces $B$ and $y$ is a local minimum of $d_x$.  The
argument from Case 1 holds unless $y$ is a vertex on an edge adjacent
to $\alpha_x$.  In such a case, we have the situation in Figure
\ref{vertint}.  The turning angle from adjacent edge containing $y$ to
$\alpha_x$ must be at least $\pi/2$.  The result follows by Lemma
\ref{mu1sep}.

\begin{figure}
  \begin{center}
    \begin{overpic}[width=3in]{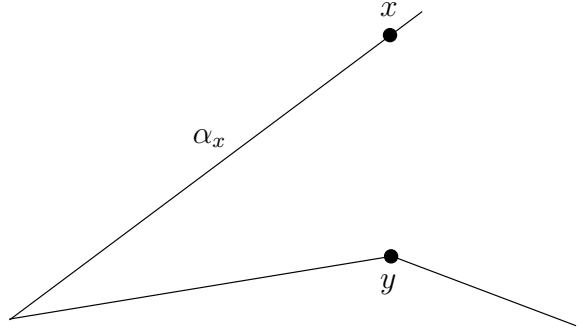}
      \put(70,70){$\alpha_x$}
      \put(141,15){$y$}
      \put(141,118){$x$}
    \end{overpic}
  \end{center}
  \caption{When $y$ is a local minimum for $d_x$ and $\alpha_x$ is
    adjacent to one of the edges containing the vertex $y$, the
    turning angle between the edges must be at least $\pi/2$.}
\label{vertint}
\end{figure}

\noindent{\bf Case 3:} $y$ is non-vertex, $x$ is a vertex\\
Let $\alpha_x^1$ and $\alpha_x^2$ be the edges adjacent to $x$.  Note
that exactly one of the adjacent edges $\alpha_x^1$ and $\alpha_x^2$
must pierce $B$ or else $(x,y)$ is doubly critical.  Without loss of
generality, assume that $\alpha_x^1$ pierces $B$.  Two things can
happen.  If $\alpha_x^1$ is adjacent to $\alpha_y$, then (as in Case
2) the turning angle must be $\pi/2$ and Lemma \ref{mu1sep} gives us
the result.  Otherwise, the argument from Case 1 holds.

\noindent{\bf Case 4:} $y$, $x$ are both vertices\\
Then $\alpha_x$ pierces $B$ and $y$ is a local minimum for $d_x$.  The
argument from Case 1 holds unless an edge containing $x$ and an edge
containing $y$ are adjacent, in which case the vertex between the two
edges must have a turning angle of at least $\pi/2$ and Lemma
\ref{mu1sep} gives us the desired result.
\end{proof}

The previous theorem says that if $dcsd(K)<MinRad(K)$, then
$scsd(K)=dcsd(K)$.  The following corollary follows immediately.

\begin{corollary}
For a polygonal knot $K$, $$MD(K) \geq \min\{2\,R(K),MinEdge(K)\}\,.$$
\end{corollary}

The following is a corollary to Lemma \ref{polytcchar} and shows that
we need only consider pairs from $DC(K)$ that are minima for the
distance function in order to characterize $R(K)$ for polygons.  Note
however, that $mdcsd(K)$ can be empty in which case we set
$mdcsd(K)=\infty$, e.g.~on a regular $n$-gon where $n$ is odd.

\begin{proposition}
Let \begin{equation*}
MDC(K)=\{(x,y)\in K\times K : x\not = y
\text{ turning point minima of }d_y \text{ and }d_x \text{ respectively}\}.
\end{equation*} 
and 
\begin{equation*}
mdcsd(K)=\min\{\Vert x-y\Vert : (x,y)\in MDC(K)\}.
\end{equation*}
Then
$$R(K) = \min\{MinRad(K),mdcsd(K)/2\}\,.$$
\end{proposition}
\begin{proof}
Note that $dcsd(K) \leq mdcsd(K)$.  If $MinRad(K) \leq dcsd(K)/2$,
then the result holds.

If $dcsd(K)/2<MinRad(K)$, then $\min_{(x,y)\in TC(K)}\Vert x-y\Vert =
dcsd(K)$.  Thus, we know that $dcsd(K)$ is realized at a pair $(x,y)$
where $y$ and $x$ are minima of $d_x$ and $d_y$ respectively, as
desired.
\end{proof}

Now we have established some properties of the polygonal thickness
radius.  Next we complete a Monte Carlo simulation to understand the
structure of the thick tube about a knot.

\section{Local Structure of Equilateral Knot Space}
A polygonal knot $K$ is determined by a sequence of $n$ points in
three-space, the vertices of $K$.  These vertices determine a sequence
of edges (straight line segments) connecting the vertices, including
an edge from the last to the first of the vertices.  It is required
that the vertices are distinct and that the edges meet only at the
common vertex with the adjacent edge sharing the vertex.  In fact,
this approach provides each knot with a first vertex and an
orientation, since the sequence determines a specific second vertex
and therefore a direction along the knot.  These knots, then, are
described by a single point in the Euclidean space of dimension $3n$
and the entire collection of such knots determines an open subset
$Geo(n)$ of this Euclidean space.  The collection of knots whose edge
lengths are all equal determines the subset of equilateral knots
$Equ(n)$ of $Geo(n)$.  Let $Edge(K)$ denote the length of the edges of
a knot $K\in Equ(n)$.  If $K$ and $K'$ are two knots in $Geo(n)$, we
define the distance between $K$ and $K'$ to be equal to the largest
standard Euclidean distance between the corresponding vertices of $K$
and $K'$.  Using this definition of distance between knots in
$Geo(n)$, we describe the neighborhood $N(K,r)$ of $K$ in $Geo(n)$ to
be those knots in $Geo(n)$ within distance $r$ of $K$.  With this
definition in mind, we observe that any perturbation of the vertices
of an equilateral knot $K$ of scale no larger than $r$ specifies a
knot in the neighborhood $N(K,r)$.  As a consequence, the analysis of
the possible knots occurring with perturbation scale less than or
equal to $r$ provides a description of the knots occurring in
$N(K,r)$.

We will need Schur's Theorem \cite{chern} for our subsequent analysis of the
necessary perturbation size required to form different types of local
knotting.

\begin{lemma}
Let $C$ and $C^*$ be two piecewise linear curves of the same length,
such that $C$ forms a simple convex plane curve with the chord
connecting its endpoints.  Let $s$ be the arclength parameter for $C$
and $C^*$.  When $s$ corresponds to a vertex of $C$ (and $C^*$), let
$angle(s)$ be the turning angle at a vertex, and $angle^*(s)$ be the
corresponding angle on $C^*$.  Let $d$ and $d^*$ be the distances
between the endpoints of $C$ and $C^*$, respectively.  Then, if
$angle^*(s)\leq angle(s)$, we have $d^*\geq d$.
\end{lemma}

Thus, given an equilateral polygon with a fixed edge length and
maximum turning angle $\theta_{max}$, we can find a lower bound for
the distance between vertices by analyzing a planar polygon with the
same edge lengths and a constant turning angle of $\theta_{max}$.

\begin{figure}
  \centering
  \includegraphics[width=3in]{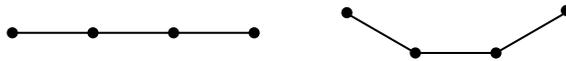}
  \caption{The perturbation size required to collect four vertices in
    a small ball changes depending on the curvature of the polygon,
    detected here by the $MinRad$ term of $R(K)$.}
  \label{fourinball}
\end{figure}

We now analyze the size of perturbation required to create local
knotting.  Note that the size of perturbation will depend on how much
the polygon is curving, namely in that more curving implies that
vertices can be closer together (see Figure \ref{fourinball}).

\begin{theorem}
If $K$ is a knot in $Equ(n)$ and the scale of perturbation $r$ is
smaller than $$\min\left\{\frac{1}{2}Edge(K)\frac{12R(K)^2 - Edge(K)^2}
{4R(K)^2+Edge(K)^2}, R(K)\right\},$$ then all knots in $N(K,r)$ are
equivalent to $K$ in $Geo(n)$.
\end{theorem}
\begin{proof}
In Section \ref{tubethmsection}, we prove the Tube Theorem (Theorem
\ref{tubethm}).  For the sake of this proof, the Tube Theorem says
that any knot change must occur within the tube of radius $R(K)$ about
the knot, what we term local knotting.  Thus, any change in knot type
must occur on an arc of $K$ whose total curvature is $< \pi$ by Lemma
\ref{tclemma}(b).

A local trefoil knot can occur at any radius beyond the value at which
four points can coalesce at a single point.  A planar gon whose radius
of curvature is equal to $R(K)$ and edge length is $Edge(K)$ has four
consecutive vertices meeting at a common point at a minimum distance
of $R_0=\frac{1}{2}Edge(K)\frac{12R(K)^2 -
Edge(K)^2}{4R(K)^2+Edge(K)^2}$.  Applying Schur's Theorem, we get that
the distance between the endpoints of a four-vertex arc is at least
this value.
\end{proof}

For an equilateral unknot with $n$ edges inscribed in a circle of
radius one, the edge length is $2\sin(\pi/n)$ and the thickness radius
is $\cos(\pi/n)$.  An analysis shows that for less than nine edges, no
local knotting is possible, i.e.~the size of the tube is smaller than
the minimum perturbation size needed to bring four consecutive
vertices to a common point.  It might be easier to think of the
thickness radius relative to the size of an edge length here.  We
include both the raw measurement of thickness radius and the thickness
radius in relation to the size of an edge.  For nine edges, this term,
0.8660 (or $1.2660\,Edge$), is smaller than the thickness radius,
0.9396 (or $1.3737\,Edge$), implying the possibility of local
knotting, i.e.~knotting occurring within the tubular neighborhood of
the knot with radius equal to the thickness radius as shown in Figure
\ref{fig1}.  In this range, only trefoil knots can arise.  This
analysis is specific to the standard regular $n$-gon, giving the
minimum perturbation required to bring together the required four
vertices as $\sin(3\pi/n)$ in order to create a local trefoil.
Applying the same method of analysis to the creation of a local
figure-eight knot where five consecutive vertices must come together,
the minimum perturbation is $\sin(4\pi/n)$.  This shows that at least
11 edges are required.  In this case the minimum perturbation is
$1.6144\,Edge$ compared to a thickness radius of $1.7028\,Edge$.  In
each of these limiting cases we should expect the proportion of the
neighborhood of the unknot consisting of trefoils and figure-eight
knots to be very very small.

To explore these and other features of the local structure, we have
undertaken a Monte Carlo sampling of the nearby knots through a random
perturbation of the vertices.  While a fuller discussion of the data
occurs later in this note, it may be helpful to introduce some aspects
of the data presently so as to give a more accurate sense of the local
structure of equilateral knot space.  For $n = 9$, one can make a
rough calculation that will give an estimate of the probability of
occurrences of trefoil summands in the perturbation.  One first
calculates the probability that four consecutive vertices coalesce in
a small ball and then multiplies this by the probability of a trefoil
in $Geo(6)$ \cite{MillettMonte}.  This crude estimate of the
probability of a trefoil summand is thereby calculated to be
$3.1974\times 10^{-16} $ or, about 3 chances in $ 10^{16} $ cases.  As
a consequence, it is not surprising that we have not encountered
trefoils in testing five million cases with nine edges.  For $n = 10$,
the probability is $8.1194\times 10^{-13}$ relatively speaking, rather
more likely, about 8 chances in $10^{13}$ cases, but still too small
to expect an observation in only five million cases.  For $n =11$,
when we first have the possibility of figure-eight knots as well as
trefoils, the theoretical lower-bound estimate of the probability of
trefoil knots increases to $6.0951\times 10^{-11}$ and, for $n = 12$,
to $1.0837\times 10^{-9}$.  The first actual observation of a trefoil
knot in our study occurred at $n = 11$ in a sample of $5\times 10^6$
cases.  This data provides a Monte Carlo estimate of the probability
of a trefoil knot at $n = 11$ within the tube of maximal thickness
radius of $ 1.2\times 10^{-6} $, rather larger than the theoretical
lower bound estimate. The first observed figure-eight knot, in
$5\times 10^6 $ trials, occurred for $n = 16$.

\begin{figure}
  \centering
  \includegraphics[width=3in]{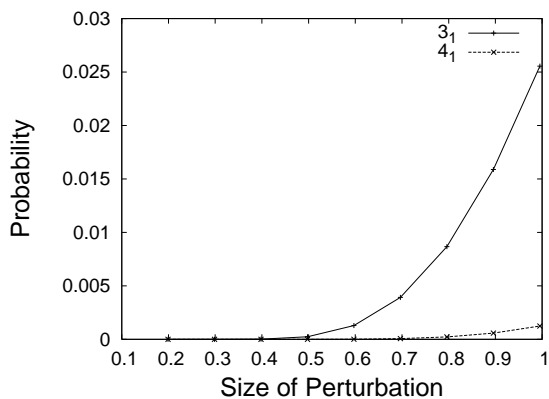}
  \caption{The proportion of perturbed $32$-gons containing exactly
    one trefoil or figure-eight knot as a function of the size of the
    perturbation.}
  \label{perturbation32.1}
\end{figure}

\begin{figure}
  \centering
  \includegraphics[width=3in]{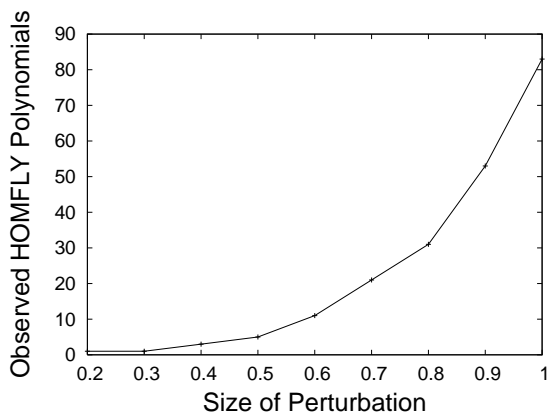}
  \caption{The number of distinct HOMFLY polynomials of knots arising
    from perturbations of the vertices of the regular 32-gon as a
    function of the size of the perturbation in a sample of size
    $2\times 10^7$.}
  \label{perturbation32}
\end{figure}

Once local knotting becomes possible, we can ask ``For a fixed number
of edges, how large of a perturbation is necessary in order that a
local trefoil might be observed?''  We did a specific analysis for the
case of regular $32$-gons.  The regular $32$-gon inscribed in the
circle of radius one has a minimum perturbation of $0.2903$
theoretically required for one or more trefoil summands, $0.3827$ for
connected sums of local trefoils and local figure-eight knots, and
$0.4714$ for more complicated summands or other satellite knots.  The
thickness radius of the regular $32$-gon is $0.9952$, approximately
one, so that a wide range of satellite knotting should be anticipated
from perturbations constrained to lie within the tube.  A Monte Carlo
exploration of the perturbations of increasing scale, with a sample
size of $2\times 10^7$, is consistent with these predictions.  The
first trefoils occur when the maximum perturbation radius is $0.3981$
and figure-eight knots and connected sums of trefoils occur when the
maximum perturbation radius is $0.4976$ (see Figure
\ref{perturbation32.1}).  The complexity in the knotting increases
rapidly with increasing perturbation size within the tube as is shown
in the graph in Figure \ref{perturbation32}.  Within the range of
perturbations remaining inside the tube, connected sums of up to three
trefoils as well as other more complex knots are observed.

\begin{theorem}
If the scale of perturbation $r$ is larger than $$\frac{1}{2}
Edge(K)\,\frac{12\,R(K)^2 - Edge(K)^2} {4R(K)^2+ Edge(K)^2}$$ but less
than $$\min\left\{4\,Edge(K)\,R(K)\,\frac{4R(K)^2 - Edge(K)^2}{
(4R(K)^2 + Edge(K)^2) ^{\frac{3}{2}}},R(K)\right\}\,,$$ then connected sums
of local trefoils may occur, but nothing more complex is possible.
\end{theorem}

\begin{proof}
Analogously, as in the previous theorem, we calculate the smallest
perturbation that allows five vertices to coalesce in a small ball,
necessary for the occurrence of a figure-eight knot.  The fact that
nothing more complex is possible is a consequence of the fact that the
perturbation limit does not allow any of the vertices adjacent to
those determining a local trefoil to meet a surrounding small ball
containing the local trefoil.  Thus, a connected sum with small
trefoils is all that can occur.
\end{proof}

\begin{theorem}
If the scale of perturbation $r$ is larger than
$$4\,Edge(K)\,R(K)\,\frac{4R(K)^2 - Edge(K)^2}{ (4R(K)^2 + Edge(K)^2)
^{\frac{3}{2}}}$$ but less than 
$$\min\left\{\frac{1}{2}Edge(K)\frac{80R(K)^4 - 40\,R(K)^2\,Edge(K)^2
+ Edge(K)^4}{(4R(K)^2+Edge(K)^2)^2},R(K)\right\},$$ then connected sum of
local trefoils and figure-eight knots may occur, but nothing more complex
is possible.
\end{theorem}

\begin{proof}
The proof of this theorem follows from the estimate of the minimal
distance needed to bring six vertices into a small ball.  It is at
this point, where we must consider bringing six or more vertices
together that more complex knots, those polygonal knots with edge
number eight, are theoretically possible \cite{calvo1, calvo3}.
\end{proof}

\section{Monte Carlo Explorations of the Local Structure of Equilateral Knot 
Space} 

As with the global character of the entirety of the knot space
$Geo(n)$, the complexity of knotting for perturbations of a knot
staying within its maximal tube increases exponentially as a function
of the number of edges. In a first study, we illustrate the structure
of local knotting by a numerical study of perturbations of the unknot.

We generated random knots within tubes about thick unknots, trefoils,
and figure-eight knots with various numbers of edges.  The host knot
is a ropelength-minimized equilateral polygon computed by TOROS
\cite{toros} using a simulated annealing algorithm.  As such, these
are not absolute polygonal ropelength-minima, but are good
approximations.  Because of the difference in the complexities of the
knots, it took many more edges to see compositions with trefoils
within the trefoil tube than the unknot tube.  Therefore, we used
different numbers of edges for the different knots.  For a given knot
and number of edges, we generated $5\times 10^6$ random knots by
perturbing each vertex by a random vector uniformly distributed in a
ball of radius $R(K)$.  The knot types were determined by computing
the HOMFLY polynomial of the perturbed knot using the Millett-Ewing
program \cite{millettewing}.  Therefore, strictly speaking, we are
computing the relative occurrence of the HOMFLY polynomial.  However,
for the range of this study, these provide an accurate representation
of the population of knots and the local incidence of these knots
parallels their occurrence in knot space.

At this point, the Millett-Ewing program can compute the HOMFLY
polynomial for knots with at most $999$ crossings.  For the
figure-eight knot with more than 500 edges, we obtained a significant
number of knots with more than 1000 crossings.  Therefore, we could
not continue the graphs beyond 500 edges for the trefoil and
figure-eight.  However, the graphs suggest that the behavior of the
relative probabilities is consistent and robust.

As is the case in studies of other knot spaces, such as $Geo(n)$,
$Equ(n)$, and the unit lattice $\mathbb{Z}^3$, the unknot is, for
small numbers of edges, the most probable knot and the probability
decays exponentially.  We fit a curve to the probabilities using the
fitting function
$$a(n-n_0)^be^{-kn-ln^2}$$ 
where $n$ denotes the number of edges and $a$, $b$, $k$, $l$, and
$n_0$ are fitting parameters.  We have tested the quality of fit
provided by a range of functions employed in various studies,
including this one, and find that this choice is appropriate for the
goals of this project \cite{MillettRawdonFitting}.  From the graphs,
once can see that the functions fit the curves well.

\begin{figure}
  \centering
  \includegraphics[width=3in]{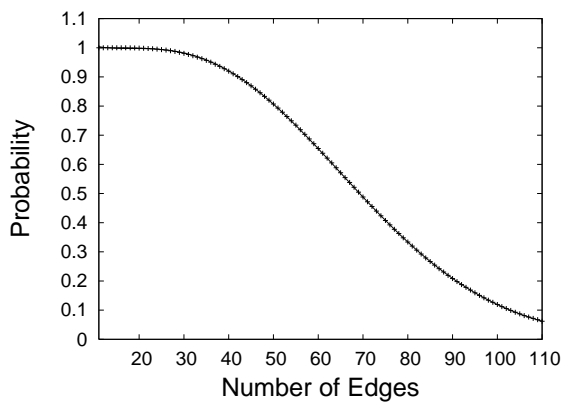}
  \caption{Probability of trivial knots versus the number of edges in
    the unknot tube.}
  \label{0.1.0.1}
\end{figure}

\begin{figure}
  \centering
  \includegraphics[width=3in]{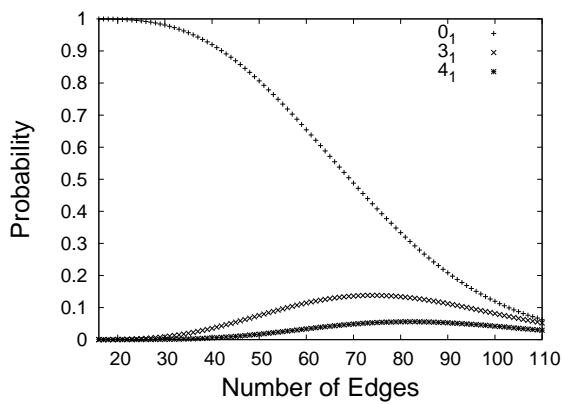}
  \caption{Probability of knotting versus the number of edges in the
    unknot tube.}
  \label{0.1.all}
\end{figure}

\begin{figure}
  \centering
  \includegraphics[width=3in]{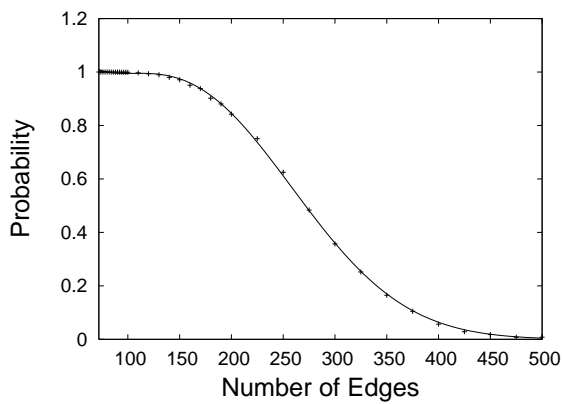}
  \caption{Probability of compositions with the trivial knot 
    versus the number of edges in the trefoil tube.}
  \label{3.1.0.1}
\end{figure}

\begin{figure}
  \centering
  \includegraphics[width=3in]{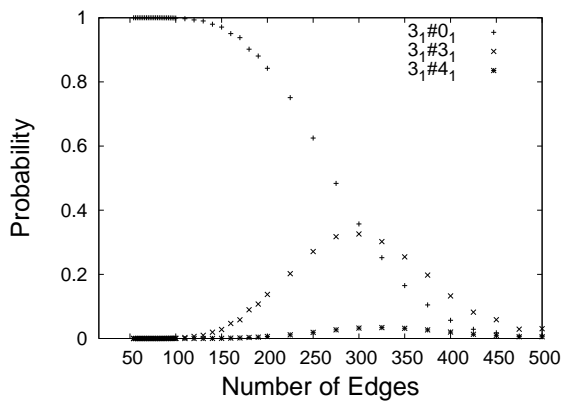}
  \caption{Probability of knotting versus the number of edges in the
    trefoil tube.}
  \label{3.1.all}
\end{figure}

\begin{figure}
  \centering
  \includegraphics[width=3in]{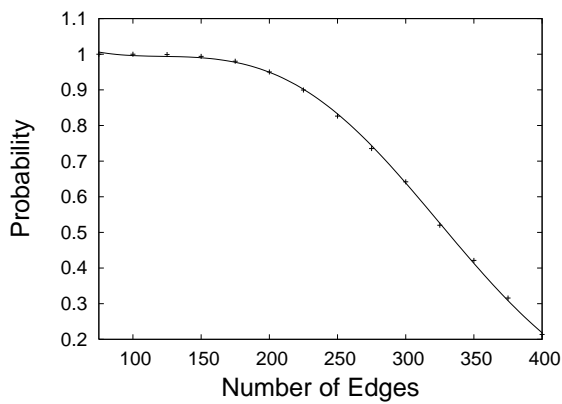}
  \caption{Probability of compositions with the trivial knot 
    versus the number of edges in the figure-eight tube.}
  \label{4.1.0.1}
\end{figure}

\begin{figure}
  \centering
  \includegraphics[width=3in]{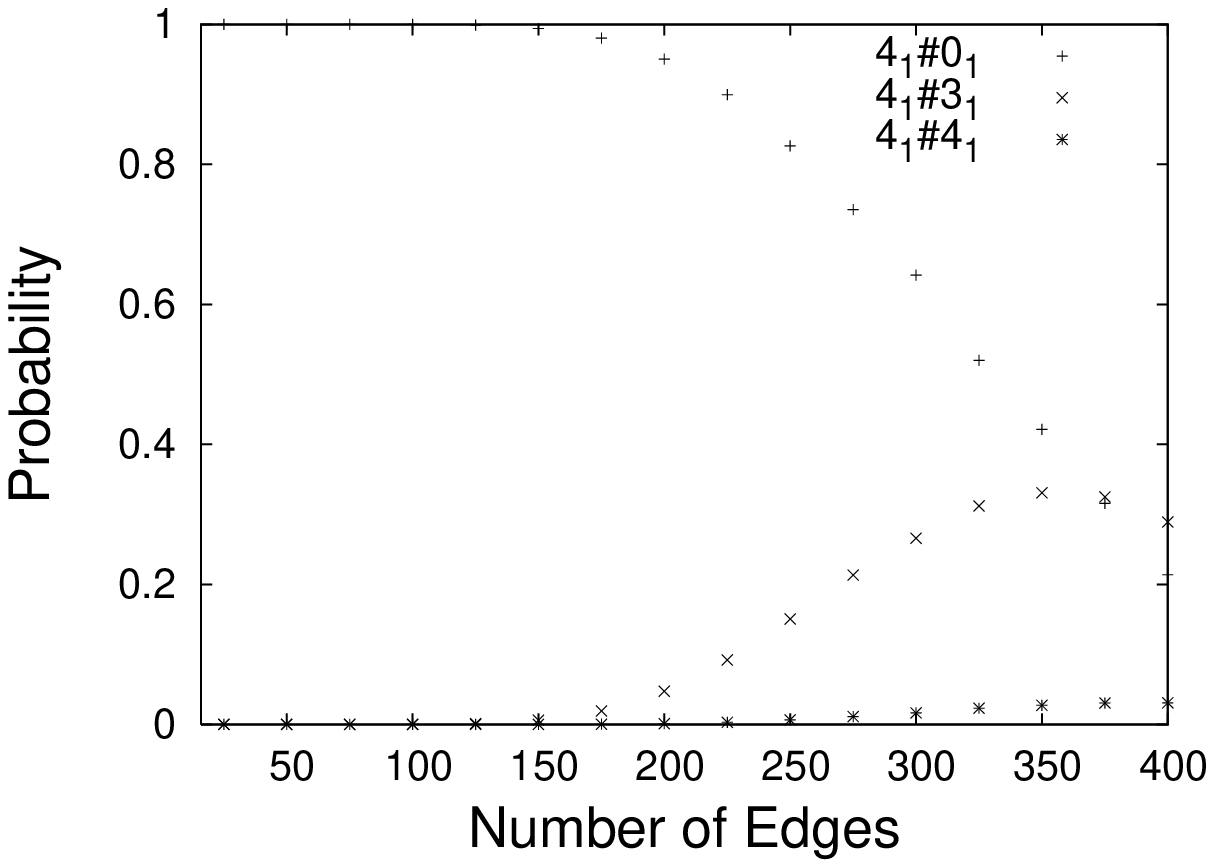}
  \caption{Probability of knotting versus the number of edges in the
    figure-eight tube.}
  \label{4.1.all}
\end{figure}

The graphs of the proportion of unknots in the thick tube about the
unknot, trefoil, and figure-eight are shown in Figs.~\ref{0.1.0.1},
\ref{3.1.0.1}, and \ref{4.1.0.1}.  The proportion of trefoil knots and
figure-eight knots, shown in Figs.~\ref{0.1.all}, \ref{3.1.all}, and
\ref{4.1.all}, demonstrate the extent of similarity of the data
regardless of the core knot.  While there are some differences, they
have remarkably similar appearances to those of the entire spaces of
geometric or equilateral knots as reported elsewhere (see
e.g.~\cite{vlfa,kamvolo,michwieg,dt1,Millettregular}).

\begin{figure}
  \centering
  \includegraphics[width=3in]{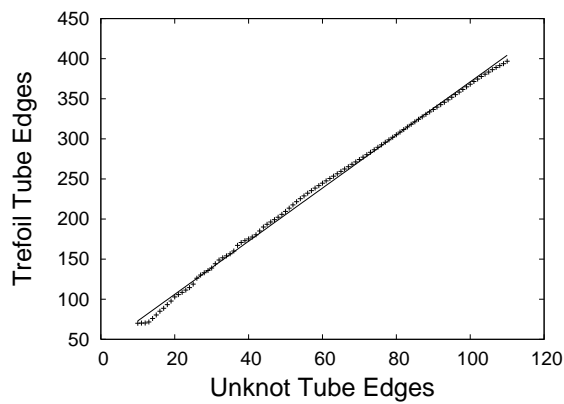}
  \caption{The correspondence between the number edges for
  perturbations within the and trefoil tubes derived by means of the
  proportion of compositions with the trivial knot.}
  \label{treftrivscale}
\end{figure}

The large scale similarity between the appearance of the local data,
that is the knot populations nearby the thick unknot, trefoil knot,
and figure-eight knot, and the equilateral and geometric knot
populations suggests that these populations may be related via a
``change of scale,'' perhaps a homeothety given by a translation and
an expansion or contraction taking $n$ to $a+b n$, where $b$
determines the change of scale.  One test of this hypothesis is to use
the monotonically decreasing data for the unknot to create a
functional correspondence between the number of edges for the
perturbations of the unknot and those of the thick trefoil.  We did
this by fitting a cubic spline to the unknot data for perturbations
within the thick trefoil knot and using a binary search to estimate
the number of edges within the unknot tube needed to achieve the same
probability.  We then fit the correspondence with a linear function.
The scale change function derived by application of this strategy to
the perturbations of the trivial knot and those of the trefoil knot is
shown in Figure \ref{treftrivscale}.  The initial behavior of the
scale change function is incompatible with a universal change of
scale.  In the range beyond this, larger than 20 edges, the unknot
data support a change of scale.  Note that for a very large number of
edges, towards the asymptotic range, a change of scale would be
expected since four consecutive vertices will be close to lying on a
straight line segment (in which case a perturbation of $1.5$ times the
edge length will be sufficient for the vertices to coalesce).  In the
present case the potential scale change function is approximately
$$ 39.8444 + 3.3127 n$$ 
with an R-squared value of $0.998$.  Looking closely at the data in
Figure \ref{treftrivscale}, one might propose a more complex structure
consisting of an initial range in which the data is sparse, a second
range until about 60 edges governed by finite size considerations, and
a subsequent range governed by asymptotic size behavior.  The
transition from finite range to asymptotic range structures remains
quite obscure and, we propose, worthy of deeper investigation.

\begin{figure}
  \centering
  \includegraphics[width=3in]{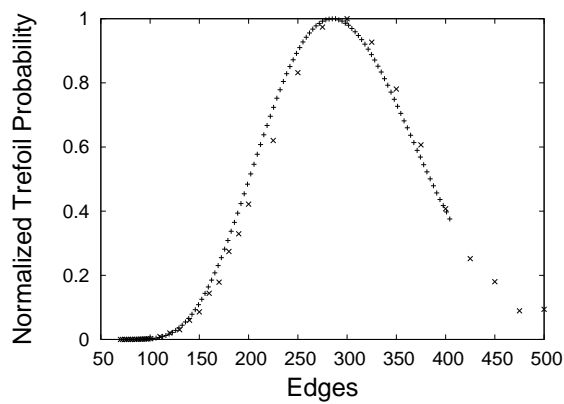}
  \caption{The trivial to trefoil perturbation change of scale applied
  to compositions of trefoils.}
  \label{trefscalev1}
\end{figure}

\begin{figure}
  \centering
  \includegraphics[width=3in]{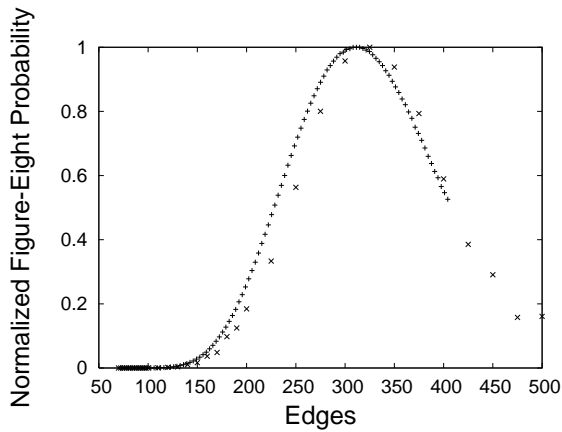}
  \caption{The trivial to trefoil perturbation change of scale applied
  to compositions of figure-eights.}
  \label{f8scalev1}
\end{figure}

A critical test of the hypothesis that a change of scale captures the
differences local and global structures of knot space is the success
or failure of the change to capture the properties of the non-trivial
knot types.  The same change of scale is applied to the trefoil knots
in Figure \ref{trefscalev1} and figure-eight knots in Figure
\ref{f8scalev1} for these perturbations.  The maximal values of the
associated knot probabilities have been rescaled to one to capture the
degree to which the shapes of these knot probability functions are
scale invariant.  While the fit is reasonable, it does not give
compelling evidence of the hypothesized scale change in as much as the
trefoil perturbation probabilities are both shifted slightly toward a
larger number of edges.  Thus, one expects that even the local knot
probabilities will be influenced by the geometric properties of the
core knot.

\section{Complexity Transitions}
\begin{figure}
  \centering
  \includegraphics[width=3in]{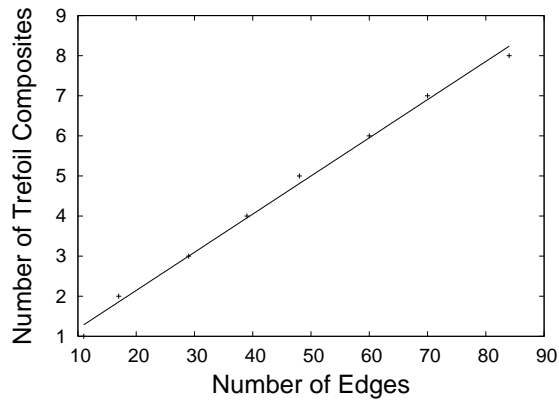}
  \caption{The observed number of trefoil summands as a function of
  the number of edges in perturbations of the unknot in samples of
  size $5\times 10^6$.}
  \label{transv1}
\end{figure}

This section concerns the increase in the complexity of observed
knotting as the number of edges increases while, roughly, keeping the
thickness constant.  For small numbers of edges the curves are so
tight that no knotting is possible, but as the number of edges
increases, knotting becomes possible.  One of the classic theorems
shows that the probability of any number of trefoil summands goes to
one and the number of edges tends to infinity.  It is in this context
that one is lead to track the number of trefoil summands that are
actually observed in the data.  One expects this number to increase as
the number of edges increases, but at what rate?  Keeping track of
this, we find that the observed occurrence is quite regular.  The
experimental data, shown in Figure \ref{transv1}, suggest that, for
the unknot, the number of trefoil summands is given by
$\lfloor\frac{n-4}{8}\rfloor$, where $\lfloor x\rfloor$ denotes the
integer part of $x$.  The fact that eight edges seems to be required
for additional trefoil summands is consistent with our earlier
theoretical estimates and the experimental data giving the first
observation of a trefoil in the perturbations of the unknot at nine
edges.  The fact that only four additional edges are theoretically
required to add another trefoil summand is offset by the exceedingly
small probability that additional knotting will occur in a random
sample.  Thus, the requirement of an additional edges is probabilistic
observation and not a theoretical consequence of the geometry.  The
extent to which these transition numbers depend upon the geometry of
the core knot is an interesting question.

\section{The Tube Theorem}
\label{tubethmsection}
Let $Tube(e_i,r)$ be the union of closed balls of radius $r$ about
each point on the edge $e_i$ and $Tube(K,r)$ be the union of
$Tube(e_i,r)$ over each edge $e_i$.  We will show that $Tube(K,r)$ is
an embedded torus when $r<R(K)$, which insures that all of the
perturbation data discussed earlier, is in fact local knotting.

Recall, in \cite{lsdr}, the thickness radius of a smooth knot was
defined as the largest radius of a non self-intersecting tube about
the knot.  The thickness radius was then characterized in terms of
local disk intersections, detected by $MinRad$, and arclength-distant
disk intersections, detected by half the double-critical self-distance
(abbreviated $dcsd$).  Polygonal thickness radius was defined
\cite{mine,meideal} to model $MinRad$ and $dcsd$ for polygons and,
thus, is based on the characterization of thickness radius for smooth
knots.  The following theorem unites polygonal thickness radius with
the origins of thickness radius in terms of tubular neighborhoods.

\begin{theorem}
If $r<R(K)$, then $Tube(K,r)$ is an embedded torus.
\label{tubethm}
\end{theorem}

\begin{figure}
  \begin{center}
    \begin{overpic}[width=3in]{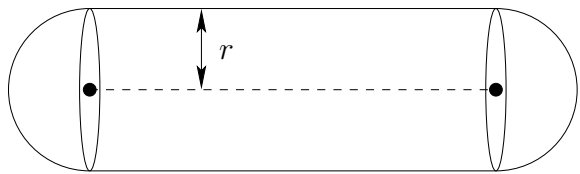}
      \put(80,43){$r$}
    \end{overpic}
  \end{center}
  \caption{The tubular neighborhood $Tube(e_i,r)$ of one edge
    (represented by the dashed line) is a cylinder with hemiball
    endcaps and resembles a pill.}
  \label{onepill}
\end{figure}

The proof of this theorem follows later in this section.  Fix an
$r<R(K)$.  We begin by analyzing the components which make up
$Tube(K,r)$.  The tube about one edge $Tube(e_i,r)$ consists of a
cylinder with two hemiballs attached at the ends and resembles a pill
(see Figure \ref{onepill}).

\begin{figure}
  \centering
  \includegraphics[height=1.4in]{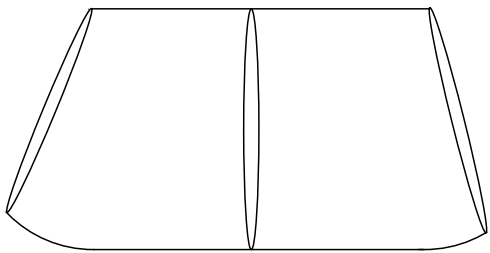}
  \caption{A cell is constructed by cutting $Tube(e_i,r)$ along angle
    bisecting planes at its two vertices.}
  \label{tiltedcell}
\end{figure}

We cut each pill along a particular angle bisecting plane at each of
its vertices to create $n$ cells and show that non-adjacent cells are
disjoint.  We construct a cell $B_i$ along $e_i$ as follows: At a
vertex $v_i$, the adjacent edges $e_{i-1}$ and $e_i$ form a plane, say
$Q_i$.  Let $P_i$ be the plane that bisects the angle at $v_i$ and
contains a normal to the plane $Q_i$.  The planes $P_i$ and $P_{i+1}$
divide space into three (in the non-generic case where two consecutive
edges are collinear) or four regions.  Let $B_i$ be the closure of
$Tube(e_i,r)$ intersected with the region which contains $e_i$.  An
example of a cell is shown in Figure \ref{tiltedcell}.

So what do these cells looks like?  We begin by dividing each cell into
several pieces.

\begin{figure}
  \begin{center}
    \vspace*{20pt}
    \begin{overpic}[height=2in]{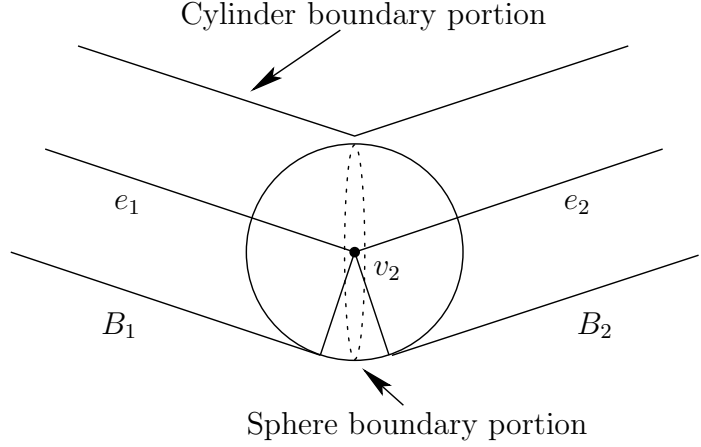}
      \put(65,145){Cylinder boundary portion}
      \put(40,75){$e_1$}
      \put(35,28){$B_1$}
      \put(138,50){$v_2$}
      \put(90,-10){Sphere boundary portion}
      \put(215,28){$B_2$}
      \put(210,75){$e_2$}
    \end{overpic}
  \end{center}
  \caption{The corner section where two cells meet.}
  \label{cornercell}
\end{figure}

Consider two adjacent cells, say $B_1$ and $B_2$ corresponding to the
edges $e_1$ and $e_2$ and vertices $v_1,v_2,$ and $v_3$ (see Figure
\ref{cornercell}).  Suppose $x$ is on the boundary of $B_1$.  Then the
minimum distance between $x$ and $e_1$ must be $r$ or $x$ is on one of
the planes that bounds $B_1$, that is $P_1$ or $P_2$.  Suppose the
former.  Then the point $y\in e_1$ minimizing the distance between $x$
and $e_1$ is either a non-vertex or a vertex point.  If $y$ is a
non-vertex point, then $x$ lies on the boundary circle of the normal
disk to $e_1$ centered at $y$ of radius $r$, i.e.~on the cylinder
boundary portion of $B_1$.  If $y$ is a vertex point, then $x$ lies on
the hemisphere of radius $r$ centered at $y$, the sphere boundary
portion of $B_1$.

We will look at the $B_1$ boundary points lying on $P_1$ and $P_2$
later.  For now, we analyze the cylinder and sphere portions of the
boundary of $B_1$.

\begin{figure}
  \begin{center}
    \includegraphics[height=1.5in]{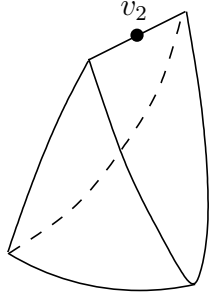}
    \put(-35,105){$v_2$}
  \end{center}
  \caption{One full spherical wedge.}
  \label{spherewedge}
\end{figure}

At the vertex $v_2$, we will show that the sphere boundary points form
a spherical wedge.  Let $H_1$ be the endcap hemiball of radius $r$ at
$v_2$ from $Tube(e_1,r)$ and $H_2$ be the endcap hemiball of radius
$r$ at $v_2$ from $Tube(e_2,r)$.  If $x$ is a vertex boundary point,
then $x$ must be an element of $H_1$ and $H_2$.  This set $H_1\cap
H_2$ is bounded by the two disks centered at $v_2$ that are normal to
$e_1$ and $e_2$ respectively (see Figure \ref{spherewedge}).  The
intersection $H_1\cap H_2$ forms a spherical wedge of interior angle
$\pi-angle(v_2)$.  The plane $P_2$ splits this wedge into two
isometric pieces, one belonging to $B_1$ and one belonging to $B_2$.
The boundary of the wedge (and the two half-wedges cut by $P_2$)
consists of points on the hemisphere and two semidisks.  Note that one
of these semidisks is connected to the cylinder of $B_1$ and the other
is connected to the other half of the spherical wedge belonging to
$B_2$.

\begin{figure}
  \begin{center}
    \begin{overpic}[height=2.5in]{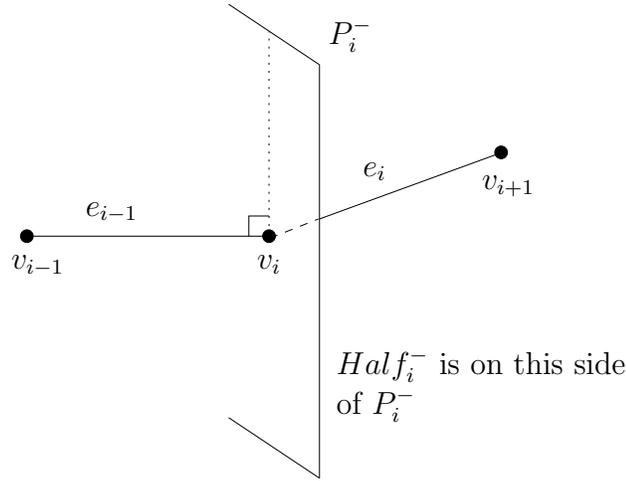}
      \put(-3,80){$v_{i-1}$}
      \put(90,80){$v_i$}
      \put(25,100){$e_{i-1}$}
      \put(120,40){$Half_i^-$ is on this side}
      \put(120,25){of $P_i^-$}
      \put(130,115){$e_i$}
      \put(175,110){$v_{i+1}$}
      \put(117,165){$P_i^-$}
    \end{overpic}
  \end{center}
  \caption{The plane $P_i^-$ is normal to $e_{i-1}$ at $v_i$.}
  \label{defnPiminus}
\end{figure}

The following lemma will play a crucial role in the proof of the Tube
Theorem.  Let $P_i^-$ and $P_i^+$ be the planes centered at $v_i$
normal to $e_{i-1}$ and $e_i$ respectively (see Figure
\ref{defnPiminus}).  Note that $P_i^-$ splits $\mathbb{R}^3$ into two
closed half-spaces, each containing $P_i^-$.  Let $Half_i^-$ be the
half-space containing all of $e_i$.  Define $Half_i^+$ similarly so
that $e_i$ is contained in $Half_i^+$.

\begin{lemma}
If $x\in Half_i^-\cap Half_i^+$, then the minimum distance from $x$ to
both $e_{i-1}$ and $e_i$ is realized at $v_i$.  In particular, if
$x\in H_{i-1}\cap H_i$, then the minimum distance from $x$ to both
$e_{i-1}$ and $e_i$ is realized at $v_i$.
\label{critvert}
\end{lemma}

\begin{proof}
It is clear that if $x\in Half_i^-$, the distance from $x$ to
$e_{i-1}$ is minimized at $v_i$.  Similarly, if $x\in Half_i^+$, the
distance from $x$ to $e_i$ is minimized at $v_i$.  Thus, if $x\in
Half_i^-\cap Half_i^+$, then the distance from $x$ to both $e_{i-1}$
and $e_i$ is minimized at $v_i$.  The second part of the theorem
follows immediately.
\end{proof}

\begin{figure}
  \vspace*{10pt}
  \begin{center}
    \begin{overpic}[height=1.5in]{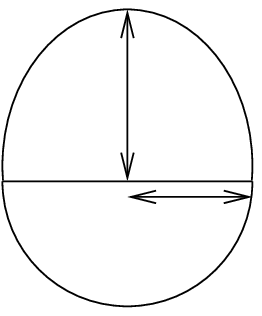}
      \put(0,112){$r\sec(angle(v)/2)$}
      \put(67,28){$r$}
    \end{overpic}
  \end{center}
  \caption{Cutting a cell along the bisecting plane at a vertex yields
    a semicircle connected to an ellipse along the diameter of the circle.}
  \label{cutcell}
\end{figure}

While the sphere boundary's two semidisk faces (just the portion of
the wedge lying in $B_1$) connect smoothly to the cylinder and other
spherical wedge, some of the cylinder boundary points connect in a
more transverse fashion.  The set of cylindrical boundary points for
$B_1$ consists of the points on the cylinder along $e_1$ minus the
portions cut off by the planes $P_1$ and $P_2$.  We concentrate on the
end with $P_2$ passing through $v_2$.  Since $r<R(K)\leq MinRad(K)$
and by the definition of $MinRad(K)$, the plane $P_2$ can only
intersect the cylinder about $e_1$ on the normal disks of the cylinder
associated to the points between the midpoint of $e_1$ and $v_2$.  In
other words, $P_1$ and $P_2$ do not intersect within $Tube(e_1,r)$.
This observation allows us to analyze this portion of $B_1$ further.
The plane $P_2$ passes through $v_2$ and bisects the angle at $v_2$.
Thus, the intersection of $P_2$ and the cylinder forms a half-ellipse
whose minor axis has length $2r$ and major axis has length
$2r\sec\frac{angle(v_2)}{2}$, which is at least $2r$.  The minor axis
is the intersection of the disk centered at $v_2$ normal to $e_1$ with
the plane determined by $e_1$ and $e_2$, and coincides with the
bounding diameter of the semidisk of the wedge.  Thus, at this end of
$B_1$, the boundary of the cylindrical portion is a semidisk with a
half-ellipse intersecting at an angle of $\frac{\pi-angle(v_2)}{2}$
(see Figure \ref{cutcell}).

The entire cell $B_1$ consists of a cylinder cut by two planes
(revealing two half-ellipses) with two spherical wedges added.  We
will show that non-consecutive cells are disjoint.  We do this by
splitting the boundary of the cell into the following components:
\begin{itemize}
\item the two half-ellipses of the cut cylinder (cylinder faces)
\item the remaining cylinder points of the cylinder boundary of the
cell (cylinder surface)
\item the two semidisk boundaries of the sphere wedge (wedge face)
\item the round spherical portion of the sphere wedge (wedge surface)
\end{itemize}

Note that the points on the cylinder surface form an open
two-manifold.  The points on the wedge surface also form an open
two-manifold.  By the discussion above, adjacent cells intersect only
along the plane bisecting their mutual vertex angle.  This
intersection is the union of a half-ellipse and semidisk (from Figure
\ref{cutcell}) joined along the minor axis and diameter respectively.

By the definition of $MinRad(K)$, normal disks of adjacent edges
centered at their respective midpoints intersect at a radius $r\geq
MinRad(K)\geq R(K)$.  Since $r<R(K)\leq MinRad(K)$, the normal disks
of radius $r$ centered at the midpoints of $e_i$ and $e_{i+1}$ do not
intersect.  Thus, the cell $B_i$ is connected and convex.  We will
paste these cells together to form $Tube(K,r)$.  We are now ready to
prove Theorem \ref{tubethm}.

\vspace*{20pt}

\noindent{\bf Theorem \ref{tubethm}.}\textit{\ If $r<R(K)$, then
$Tube(K,r)$ is an embedded torus.}

\begin{proof}
Clearly, $Tube(K,r)=\cup_i B_i$.  We will show that non-consecutive
cells are disjoint.  Let $r$ be the minimal radius such that the
intersection of a pair of non-consecutive cells is non-empty.  We show
that either $r\geq R(K)$ or $r$ is not minimal, to obtain a
contradiction.  Suppose $B_i$ and $B_j$ ($i,j$ non-consecutive)
intersect.  Since $r$ is minimal, this intersection must occur on the
boundary of $B_i$ and $B_j$.  Let $x\in B_i\cap B_j$.

The remainder of this proof is split into several cases, each one
covering a situation where $x$ is in a different part of the boundary
of $B_i$ and $B_j$.  First, we cover the intersections between
surfaces.  Note that all of these points lie on the boundary of the
cylinder or spherical wedge and are distance exactly $r$ from $K$.

\noindent{\bf Case:} Cylinder surface to cylinder surface 

If the tangent planes of the cylinders at $x$ are not identical, then
one of the cylinders pierces the other cylinder, that is, the
interiors of the three-dimensional filled cylinders intersect.  Thus,
there is an intersection between the cylinders at a radius smaller
than $r$.  This contradicts that assumption that $r$ is minimal.

If the tangent planes coincide, then there exists $x_i\in e_i$ and
$x_j\in e_j$ such that $x$, $x_i$, and $x_j$ are collinear and
$e_i\perp \overline{x_ix_j}\perp e_j$.  Thus, $(x_i,x_j)$ is a doubly
critical pair.  This implies $\Vert x_i-x_j\Vert \geq 2R(K)$ and
$r\geq R(K)$.

\noindent{\bf Case:} Cylinder surface to wedge surface

As above, if the tangent planes do not coincide, then the cylinder
pierces the wedge.  This contradicts the fact that $r$ is minimal.

So suppose the tangent planes coincide.  Without loss of generality,
we can assume that the wedge is a part of $B_i$ and the cylinder is a
part of $B_j$.  Then there exists a vertex $v\in e_i$ (whichever of
$v_i$ or $v_{i+1}$ determine the spherical wedge) and a point $x_j\in
e_j$ such that $\overline{vx_j}\perp e_j$.  By
Lemma \ref{critvert}, the pair $(v,x_j)$ is actually doubly critical.
This implies $\Vert v-x_j\Vert\geq 2 R(K)$ and $r\geq R(K)$.

\noindent{\bf Case:} Wedge surface to wedge surface

If the tangent planes do not coincide, then the wedges intersect at a
smaller radius and $r$ is not minimal.

Suppose the tangent planes coincide.  Let $u$ and $v$ be the vertices
corresponding to the wedges.  By Lemma \ref{critvert} (applied twice), 
$(u,v)$ forms a
doubly critical pair and we have $r\geq R(K)$.

\noindent{\bf Case:} Wedge face to wedge surface and cylinder surface

The semidisk face of the half-wedge intersects smoothly with the
adjacent edge's half-wedge.  Thus, any intersection on this semicircle
boundary is the same as a wedge surface intersection.  Any
intersection of the remaining portion of the disk must intersect the
interior of one of the two adjacent half-wedges.  This would
contradict that $r$ is minimal.

\noindent{\bf Discussion of the cases:} Cylinder face to anything

The intersection of the half-ellipse and the semidisk, a line segment
which is the minor axis and bounding diameter respectively, can be
considered as a part of the semidisk and are covered by the argument
above.  We restrict our attention to the remaining half-ellipse
points.

We will show that nothing can intersect the points on the half-ellipse
without contradicting the assumption that $r$ is minimal.  The
cylinder face corresponds to two adjacent cut cylinders, say $B_i$ and
$B_{i+1}$, so we will consider the region formed by the union of their
cut cylinders (which is bounded, in part, by the boundary of the
common half-ellipse).  If the turning angle at the vertex is 0, then
the adjacent edges are collinear and the union of the two cells is
cylindrical in a neighborhood of the half-ellipse (which would be a
semicircle in this case).  So suppose this angle is not 0.  Then the
tangent planes at any half-ellipse boundary point do not coincide.  If
$x$ is a half-ellipse point, then any disk centered at $x$ lies
partially within either $B_i$ or $B_{i+1}$ (or both).  In other words,
there is a crease between the adjacent cells.  We continue with the
four remaining cases.

\noindent{\bf Case:} Cylinder face to cylinder or wedge surface

Let $D$ be a small disk centered at $x$ lying on the tangent plane of
the cylinder or wedge surface at $x$.  Then $D$ intersects the
interior of $B_i$ or $B_{i+1}$ at a radius smaller than $r$, which is
a contradiction.

\noindent{\bf Case:} Cylinder face to wedge face

Again, the faces of the wedge can attach smoothly to the cylinder and
wedge surface.  Thus, this case is covered by the previous argument.

\noindent{\bf Case:} Cylinder face to cylinder face

Let $D$ be a disk centered at $x$ lying entirely within $B_i$ (there
are many of these).  By the observation above, this disk must
intersect the interior of $B_j$.  Thus, the interior of $B_i$
intersects the interior of $B_j$, which implies that $r$ is not
minimal.
\end{proof}

The previous result is sharp in the sense that for any polygonal knot
in which $R(K)=dcsd(K)/2\leq MinRad(K)$, any slightly larger
neighborhood will necessarily no longer be an embedded torus.  While
this result has not been shown for polygons, Durumeric \cite{ozlocal}
has shown that ropelength-minimized smooth knots necessarily have
$dcsd(K)/2\leq MinRad(K)$.  Furthermore, our construction yields a
canonical strong deformation retract from the tube to the knot.  This
can be obtained by ``fanning'' over the points in the cut portion of
the cells and sending the remainder of the points (i.e.~points which
lie on a normal disk where all of the normal disk lies inside a
particular cell) to the center of the disk.

\begin{figure}
  \centering
  \includegraphics[height=1.5in]{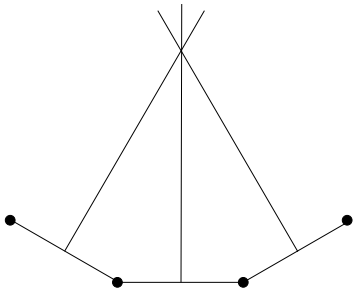}
  \caption{When $r$ is larger than $MinRad(K)$, the tube extends beyond
    the intersection of the perpendicular bisectors of the edges.}
  \label{rtoobig}
\end{figure}

However, if the thickness radius is controlled strictly by the
$MinRad$ term, a radius $r$ slightly larger than $R(K)$ can still be
an embedded torus.  If $r>MinRad(K)$, then a pair of consecutive
planes, say $P_i$ and $P_{i+1}$, will meet in the interior of
$Tube(K,r)$ (see Figure \ref{rtoobig}).  While $Tube(K,r)$ could still
be a embedded torus, we lose the canonical strong deformation retract
from this torus to the knot.  The smooth thickness radius has this
same behavior.

Furthermore, it is possible that for some $r>R(K)$, the tube
$Tube(K,r)$ is an embedded torus.  Imagine a small, tight trefoil in
what would otherwise be a large circle.  When $r$ is sufficiently
large to ``fill-in'' the volume around the trefoil, there is an
interval of radii for which the solid will be an embedded torus.  The
radii at which such phase changes occur may be of independent
interest.  In particular, these critical radii determine scales at
which different satellite structures of local knotting can appear.

\section{Conclusions}

The theoretical and computational explorations of knots near thick
equilateral trivial, trefoil, and figure-eight knots in the space of
geometric knots provide a clear picture of the
growth of complexity.  While small perturbations give topologically
equivalent knots, with increasing size of perturbation the complexity
of knotting encountered grows rapidly all the while remaining within
the tubular neighborhood of the knot.  We have shown, in Theorem
\ref{tubethm}, that this neighborhood can be understood as the union
of balls centered at points of the core knot and have established
relationships between the length of the edges and the thickness radius
of the knot that gives structure to the increase of complexity with
increasing size of perturbation.  These knots are all satellites of
the core knot and are most likely, as is demonstrated by our
simulation data, to be connected sums of the core with local knots in
the tube.

A numerical study of the growth of knotting about the thick knots
shows significant qualitative similarity that might be an indication
of a change of scale between the neighborhoods.  Testing this
hypothesis with the cases of the trivial and trefoil tube perturbations, we see
evidence of this possibility but, also, a distinct shift in structure
that indicates that this may not be true due to the finite range
behavior of the knotting distributions.

The numerical study confirms and quantifies the structure of local
summands by testing the growth of trefoil summands in the
perturbations of the trivial knot.  With each
additional eight edges one is likely to see the possibility of an
additional trefoil sum in a sample of size $5\times 10^6$.

There are several interesting questions that are suggested by the
results of this project.  First, the extent to which the geometry of
the core knot influences the complexity of the knotting arising from
small perturbations is suggested in our work but, we believe, there
may be stronger connections to be uncovered even in the case where the
knot structure is relatively homogeneous.  Secondly, our work on
perturbations describes knotting near a thick knot in the space of
geometric knots but does not have any direct implications on the
neighborhood of the knot in the space of equilateral knots.  Such a
study would be interesting and important as we would expect
the limitation to equilateral knots to promote a much more rigid
structure than that reflected in the free perturbations.  Third, as we
have discussed elsewhere \cite{MillettRawdonFitting}, the knot
data in this study lies within the ``finite range,'' arguably before the
domination of the exponential decay that defines the asymptotic range
of knotting.  Already, in this finite range, one finds great
complexity in the knotting structure.  Explorations that give high
quality data in the asymptotic range may, one imagines, give a
qualitatively different picture.  Fourth, in this work we have
investigated the structure of ``local knotting,'' the summands which
occur within the tube surrounding the core knot.  We have not yet
observed non-summand satellites that are theoretically possible but
unlikely.  What happens, we wonder, if one moves into the global
knotting regime by allowing perturbations that extend beyond the
frontier of the tube.

\section*{Acknowledgments}
Piatek and Rawdon were supported by NSF Grant \#DMS0311010.

\bibliographystyle{plain} 
\bibliography{wonderbibv12}

\end{document}